\documentclass[12pt]{article}
\usepackage{a4,verbatim,amsmath,amssymb}
\usepackage[isolatin]{inputenc}
\usepackage{epsfig}
\input math1.def
\usepackage{amsthm}

\newtheorem{theor}{Theorem}[section]

\newtheorem{lemma}[theor]{Lemma}

\numberwithin{equation}{section}
\def\proof{\goodbreak\noindent{\sc Proof. }\nobreak}

\def\endproof{\par\nobreak\hbox to \hsize{\hfil\vrule width 5pt height
5pt}\goodbreak\vskip 3pt}

\font\textmsbm= msbm10 scaled 1200
\font\scriptmsbm= msbm7 scaled 1200
\font\scriptscriptmsbm= msbm5 scaled 1200

\font\textmsbm= msbm10 scaled 1200
\font\scriptmsbm= msbm7 scaled 1200
\font\scriptscriptmsbm= msbm5 scaled 1200
\newfam\tmwfama
\textfont\tmwfama\textmsbm
\scriptfont\tmwfama\scriptmsbm
\scriptscriptfont\tmwfama\scriptscriptmsbm
\def\bbb{\fam\tmwfama\textmsbm}
\def\bT{{\bbb T}}

\def\bC{{\bbb C}}
\def\bN{{\bbb N}}

\def\bD{{\bbb D}}
\def\bE{{\bbb E}}
\def\bF{{\bbb F}}

\def\bN{{\bbb N}}

\def\bP{{\bbb P}}

\def\bT{{\bbb T}}

\def\a{\alpha}
\def\b{\beta}

\def\d{\delta}
\def\e{\epsilon}

\def\r{\rho}

\def\s{\sigma}

\def\t{\tau}

\def\vp{\varphi}
\def\o{\omega}

\def\O{\Omega}

\def\O{\Omega}

\def\sb{\subset}      \def\sbe{\subseteq}

\def\bydef{\,\lower-.1ex\hbox{\rm :}\!=}

\def\cF{{\cal F}}

\def\cH{{\cal H}}

\begin{document}
\title{A Decomposition for  Hardy Martingales  III}
     \author{Paul F. X. M\"uller\thanks{Supported
by the Austrian Science foundation (FWF) Pr.Nr. FWFP23987000. 
Participant of the Hausdorff Trimester on  Harmonic Analysis and PDE in Bonn, 2014  }  }
\date{April $8^{th}$, 2015}
\maketitle
\begin{abstract} 
We prove  Davis  
decompositions  for vector valued   Hardy martingales and illustrate their use.
This paper continues \cite{pfxm12} and \cite{pfxmpart2} on Davis and Garsia  
Inequalities.
 
\end{abstract}
\tableofcontents
\newpage

\section{Introduction}\label{intro}

The book by  A. Pelczynski \cite{pelc},
 ``Banach Spaces of analytic functions and absolutely summing operators, 
(1977)''
 contains  -inter alia- the following problems:  
\begin{enumerate}
\item 
 Does $H^1 $ have an unconditional basis?
\item
 Does there exist a subspace of $L^1 / H^1 $ isomorphic to $ L^1 $?
\item
 Does  $L^1 / H^1 $ have cotype 2?
\item
 Are the spaces $ A(\bT ^n )$ and $ A(\bT ^m )$ not isomorphic when $ n  \ne $  $ m $ ? 
\end{enumerate}
It is well known that the solutions to 
those four  problems were obtained by B. Maurey \cite{m1} and 
J. Bourgain \cite{b1, b2, bourgain-acta} respectively.
 A common feature of the proofs by Maurey and  Bourgain 
is the systematic use   of certain complex analytic martingales. Those were  studied in detail by  
D.J.H. Garling \cite{ gar2, garlingyudin} who 
coined the term {\it Hardy martingales}. 
\paragraph{Scalar valued Hardy martingales} were developed from different
viewpoints  
by B. Maurey \cite{m1} and  by J. Bourgain \cite{b1} to obtain 
the isomorphisms 
that gave the solution of the  first two problems.   

Motivated by  \cite{b1} we showed recently \cite{pfxm12}, 
that any scalar valued Hardy martingale
$F$ may be decomposed  as $F = G+ H $ into the sum of two Hardy martingales 
so that 
\begin{equation}\label{int1}
 | \Delta G_k | \le C | F_{k-1} |
\quad\text{and}\quad 
 \bE \sum |  \Delta B_k | \le C \| F \| _{L^1} ,
\end{equation}  
where $ \Delta B  _k  =  B_k  - B_{k-1} $ and $ \Delta G  _k  =  G_k  - G_{k-1}. $
We used  a non-linear telescoping trick  to derive from \eqref{int1}
the Davis and 
Garsia inequalities for scalar valued Hardy martingales 
\begin{equation}\label{int2}
  \bE ( \sum \bE_{k-1}  | \Delta G_k |^2 )^{1/2}  +   
\bE \sum |  \Delta B_k | \le C \| F \| _{L^1} .  \end{equation} 
The estimates \eqref{int1} and \eqref{int2} are specific for Hardy martingales 
and cease to hold true in general. In \cite{pfxmpart2} we determined 
 the extent to which  \eqref{int1} and \eqref{int2} are stable under dyadic perturbations, and described 
the role of the perturbed estimates in the proof that 
 $ L^1 $ embeds into  $L^1 / H^1 .$   
\paragraph{Vector valued Hardy martingales} were crucial  in the solution of problems 
3. and 4.
The martingale inequalities  that  gave rise to the cotype 2 property 
of the quotient space   $L^1 / H^1 $,  and the  isomorphic invariant 
that distinguishes between the polydisk algebras in  
different dimensions, 
are  expressed in terms of vector valued Hardy martingales. See \cite{bourgain-acta, b2, bourgaindavis}. 
Bourgain's isomorphic 
invariant \cite{b2} quantifies the fact that Hardy martingales ranging in 
the dual spaces   $ A^*(\bT^n )$ respectively  $ A^*(\bT ^m )$ behave 
significantly different when $ m \ne n . $ The vector valued 
Riesz product studied by G. Pisier (see \cite{dgt}) gave rise to  Hardy martingales 
with values in   $L^1 / H^1 $ that intertwine the cotype 2 properties 
of   $L^1 / H^1 $, and Bourgain's isomorphic invariants in \cite{b2}.
It also  played an important role for  the work of 
{ W. Davis, D. J. H. Garling, N. Tomczak-Jaegermann} \cite{dgt} on Hardy martingale cotype and complex uniformly convex 
renormings of Banach spaces.  
\paragraph{In the present  paper}  we study decompositions for vector valued 
Hardy martingales.
Our point of reference is the following theorem of B. Davis \cite{davis}. 
If an $X$ valued martingale $ F = ( F_k ) $ satisfies 
$$ \bE (\sup_{k \in \bN } \| F_k\|_X ) <\infty,$$
then there exist martingales  $ G = ( G_k ) $ and  $B = ( B_k ) $
with $ F _ k = G _k + B _ k, $ $ k \in \bN  $ so that 
\begin{equation}\label{in1}
 \| \Delta G _k \|_X  \le C \max _{m \le k-1} \| F_m \|_X  \quad
\text{and } \quad 
\bE \sum   \| \Delta B  _k \|_X \le C  \bE (\sup_{k \in \bN } \| F_k\|_X ),  
\end{equation}
The vector valued decomposition theorem of B. Davis  
and our  previous work on scalar valued Hardy martingales,  \cite{pfxm12} and \cite{pfxmpart2}, 
gave rise to  the following questions:
\begin{enumerate} 
\item  Is it still possible to prove  this  decomposition 
under the additional constraint    that $ F$ and  $ G , B $ 
are vector valued Hardy martingales? 

It is easy to see, and well known, that  the original 
proof by Davies \cite{davis, sia} does not preserve the class of Hardy martingales.
In this paper  we therefore use  a new decomposition that respects the 
condition of  analyticity  and simultaneously 
yields  the estimates \eqref{in1}. See Theorem  \ref{dc1}. The construction is  based on Brownian motion stopping times and 
Doob's martingale projection operator. 

We apply  the  decomposition theorem \ref{dc1}
to prove an extrapolation result   for operators acting on 
Hardy martingales. Theorem \ref{dc2}.
\item 
Is it possible to further exploit that $F$ is taken in the 
class of Hardy martingale and obtain   an improved   decomposition
with  estimates that go beyond  \eqref{in1}? 

In response to this question   in Theorem \ref{dc4}  we obtain 
a decomposition   of $F$  into Hardy martingales $F = G + B ,$ 
for which we prove
the following estimates 
\begin{equation} \label{in2}
 \| \Delta G _k \|_X  \le C  \| F_{k-1} \|_X  \quad
\text{and } \quad 
\bE \sum   \| \Delta B  _k \|_X \le C  \bE  \| F\|_X .  
 \end{equation}
The splitting itself is done again by  Brownian motion, stopping times and Doob's projection; 
the verification of  \eqref{in2} relies on  Havin's lemma and outer functions.

The estimates  \eqref{in2}  and \eqref{in1} hold true for any complex Banach space;
thus Hardy martingales are to general martingales as \eqref{in2} is to \eqref{in1}.

With the decomposition estimates \eqref{in2}  and  hypothesis ``$\cH(q)$''  we obtain 
 further inequalities for vector valued   for Hardy martingales.
Let $ q \ge 2 . $ A   Banach space $X$ satisfies the hypothesis $\cH(q)$  if for each 
$ M \ge 1 $ there exists $ \d = \d ( M ) > 0 $ such that for any  $ x \in X $ with $ \| x \| = 1 $ and
$ g \in H^\infty _ 0 ( \bT , X )  $ with $ \| g \| _\infty  \le M , $ 
\begin{equation}
\label{7-4-15-5} \int_\bT \| z + g \|_X  dm  \ge  (  1  + \d (\int_\bT  \|  g \|_X ^ q   dm )^q )^{1/q}
  . \end{equation} 
Theorem \ref{6-4-15-26} asserts that if the  Banach space $X$ satisfies  $\cH(q)$  
then any   $X$-valued  Hardy  martingale $F$ has a decomposition into Hardy martingales as $F = G + B $
such that
\begin{equation*}
\bE (  \sum _{k= 1 }^\infty  ( \bE_{k-1}\|\Delta G_k \|_X  ^q ))^{1/q}    +   \bE ( \sum _{k= 1 }^\infty  \|\Delta B_k \|_X ) 
 \le A_q \bE( \|F \|_X) .  \end{equation*} 
If we replace \eqref{7-4-15-5} by the weaker  hypothesis
\begin{equation}
\label{7-4-15-45}  
 \int_\bT \| z + g \|_X  dm  \ge  (  1  + \d (\int_\bT  \|  g \|_X  dm )^q )^{1/q}
  , \end{equation} 
then  we are able to prove that the decomposition estimates \eqref{in2}  yield 
\begin{equation*}
  \bE (  \sum _{k= 1 }^\infty  ( \bE_{k-1}\|\Delta G_k \|_X  )^q )^{1/q}    +   \bE ( \sum _{k= 1 }^\infty  \|\Delta B_k \|_X )
\le A_q \bE( \|F \|_X) .  \end{equation*} 
We note in passing that for scalar valued analytic functions, when   $ X = \bC$,  the conditions     \eqref{7-4-15-5} and  \eqref{7-4-15-45} hold true with $ q = 2  .$   See    \cite{b1, pfxm12}.

\item Illustrating the use of Brownian Motion we give a simple proof of the fact that any Hardy martingale can be embedded -as a subsequence-  
into another  Hardy martingale with small and previsible increments.
Theorem \ref{dc7} should probably be regarded as a weak version of Q. Xu's  embedding theorem referred to by  Garling \cite{gar2}. 
See also the construction of  G. Edgar \cite{edgar1, edgar2}.  
\end{enumerate}

\section{Preliminaries}

\paragraph{Hardy spaces. }
Let  $ X $ be a complex Banach space. For  $ 1 \le p \le \infty  $ we denote by  $L_0^p(\bT , X ) ,$ 
 the Bochner space of of $X$ valued $p- $  integrable,
 functions with vanishing mean. Here  $\bT=  \{ e^{i\theta} : \theta \in
[0, 2\pi [ \} $ is   the torus 
equipped with the  normalized angular measure.
We define   $H^p_0 (\bT, X  ) \sb  L^p_0 (\bT, X  ) $ to   consist of those 
functions for which the harmonic extension to the
unit disk is analytic.  
See  \cite{p1} ,\cite{edgar1}, \cite{gar2}.

\paragraph{Martingales on  $\bT^\bN$. }
 Let $\bT^\bN $
be the   countable torus product
 equipped with its normalized  Haar measure $\bP .$
We enote by $\cF_n $
the  sigma-algebra on $\bT^\bN$ generated 
$ \{(A_1, \dots, A_n , \bT^\bN )\},$
where $A_i,\, i \le n $ are measurable subsets 
of $\bT . $ 
Let   $F = (F_n) $  be a sequence  in the Bochner space  $L^1(\bT^\bN, X )-$ 
so that $ F_n$ is $\cF_n$ measurable.  It  is an   $(\cF_n) $ martingale
if conditioned on  $\cF_{n-1}$  the  difference $\Delta F_n =
F_n- F_{n-1}$  defines an element in  $L_0^1(\bT, X ) .$
Doob's maximal function estimate states that
\begin{equation}\label{linz1}
 (\bE  \sup _{k \in \bN }  \| F_k \|_X^p )^{1/p }   \le \frac{p}{p-1}  ( \sup _{k \in \bN }   \bE  \| F_k \|_X ^p )^{1/p }  , \quad \quad  1 < p \le \infty ,
\end{equation}
for every  $X$ valued   $(\cF_n) $ martingale. ( See e.g. \cite{durr}. )


Assume now that  $F = (F_k) $ is an  $X$ valued   $(\cF_n) $ martingale.
It is called a {\it Hardy martingale} if conditioned to  $\cF_{n-1} ,$ 
the martingale difference  
$ \Delta F _n  = F_n - F_{n-1}  $
defines an element in $H^1_0 (\bT, X  ).  $  See Garling  \cite{gar2}, Pisier \cite{p1}.

\paragraph{Brownian motion.}
 Let $ \O $ denote the Wiener space.  
We let $ \{ z_t : t > 0 \} $ denote complex Brownian motion started at 
$ 0 \in \bD , $  let  $ \{\cF_t  : t > 0 \}    $ denote its  
associated continuous 
filtration, and define the stopping time $ \t $ to denote the 
first time when  Brownian motion $ \{ z_t : t > 0 \} $ hits 
the boundary of the unit disk, thus  
$$ \t = \inf \{ t > 0 : | z_t | > 1 \} . $$
We recall that for any $ f \in  H^1 (\bT , X) $ and 
$ 0 < \a \le 1, $ 
$(  \|f( z_{t\wedge\tau } )  \|_X  ^\a )$ is a submartingale, and that Garling's inequality  \cite{gar2} asserts that 
 $$  \bE (  \sup_{t < \t } \| f( z_t ) \|_X ) \le e \sup _{t < \t }    \bE ( \| f( z_t ) \|_X ) , $$
where the integration  is taken over the Wiener space  $ \O .$
We recall   {\it Doob's projection operator}  $N : L^p (\O , X) \to    L^p (\bT , X  ) $ 
 acting, 
by  conditional expectation,  on random variables defined on
Wiener space $\O ,$  
$$ N f( z ) = \bE ( f |   z_{\t}  =  z ) , \quad \quad  z \in \bT .
  $$
We use Doob's martingale projection to generate analytic functions in  $H^\infty (\bT , X  )  $
by the following stopping time procedure.
For $ f \in  H^1 (\bT , X  ) , $  $\lambda > 0$  put  
$$  \rho   = \inf \{ t < \t : \|f(   z_{t}  ) \|_X >  \lambda \} , \quad\quad 
 R =  f(   z_{\rho }  ),
 \quad\quad g = N(R) .  $$
Then $ g $ is analytic and uniformly bounded by $  \lambda . $ Precisely, 
\begin{equation}\label{him-1}
 \| g\|_X \le \lambda, \quad\quad g \in  H^\infty (\bT , X  ) .\end{equation}
See \cite{var1} for the original argument based on duality, and \cite{jopfxm} 
for an
alternative proof, based on Ito calculus.
\paragraph{Maximal function estimates for Hardy martingales \cite{gar2}.}  
Let $X$ be a Banach space and let $ F = ( F_k ) $ be an integrable $X$ valued 
Hardy martingale. For any $ 0 < \a \le 1  , $  
$( \| F_k \| _X ^\a)  $ is a non-negative submartingale,
$$ \|  F_{k-1} \| _X ^\a \le \bE_{k-1} ( \|  F_k \| _X ^\a ) , $$
and
\begin{equation}\label{him0}
 \bE ( \sup _{k \in \bN }  \| F_k \|_X )  \le e  \sup _{k \in \bN }   \bE ( \| F_k \|_X ) . 
\end{equation}
Moreover for any $ k \in \bN ,$  Garling's  theorem \cite{gar2} yields that 
 the Brownian maximal function 
 $$ F_k ^* ( x , \o ) = \max 
\left\{ \max _{ m \le k-1 } \| F_{m } ( x ) \| _X , \sup _{t < \t } 
\| F_{k} ( x , z_t ( \o ) ) \| _X  \right \}, \quad  x \in  \bT^{k-1} $$   
is integrable over $ \Sigma =  \bT^{k-1} \times \O  $  and   
\begin{equation}\label{him1} \bE _\Sigma (  F_k ^* ) \le C \bE ( \| F_k \| _X ) .
\end{equation}

\section{Vector Valued Hardy Martingale Decompositons}

\subsection{The classical Davis decomposition}

Here we present  martingale decompositions that preserve the class of vector
valued Hardy martingales.   
We split such an  $F$  as  
$F = G + B $ where $G $ is a vector valued Hardy martingale with predictable 
increments, and where the martingale differences   of $B$ are  absolutely summing. 
The proof  combines 
Davis's original  idea and 
maximal function estimates \eqref{him0},  \eqref{him1} and  the fact that 
Doob's  projection  $N$ preserves analyticity \eqref{him-1}.  

\begin{theor}\label{dc1} Let $X$ be a Banach  space. 
Any  $X$ valued  Hardy martingale  $F = (F_k)_{k = 1 }^n $
can be decomposed into the  sum of  Hardy martingales $ F = G + B $ 
such that
$$ \|\Delta G_k \|_X \le 2 \max_ {m \le k-1} \|F_m\|_X , $$
and
$$ \bE ( \sum _{k= 1 }^n  \|\Delta B_k \|_X ) \le C \bE( \|F \|_X) . $$
\end{theor} 
\proof
Fix $ k \in \bN   $  and  condition 
 to $ \cF_{k-1} . $ That is we fix $ x \in \bT^{k-1} , v \in \bT $ and put
$$ f( v ) =   \Delta F_{k}  ( x , v ) , \quad \quad \lambda =  \max_ {m \le k-1} \|F_m ( x ) \|_X . $$
Define 
$$  \rho   = \inf \{ t < \t : \|f(   z_{t}  ) \|_X > 2 \lambda \} , \quad\quad 
A = \{  \rho  <  \t \} .  $$
Now put 
$$ R_k =  f(   z_{\rho }  )  , \quad\quad  S_k = 
     f(   z_{\rho }  ) -  f(   z_{\t }  ) . $$

We next analyse the properties of  $ R_k $, $G_k .$
For $ \o \in A ,$ 
$$ F_k ^* (x,  \o) \le 4(  F_k ^* (x,  \o) -  F_{k-1} ^* (x,  \o)) .$$
and by definition  $ S_k $ is supported on $A,$ hence  
$ \| S_k \|_X 
\le 2 F^* _k \le 
8 (  F_k ^* -   F_{k-1} ^*  ) , $
and 
\begin{equation}\label{him2} \sum_{k=1} ^n  \| S_k \|_X  \le 8   F_n ^* . 
\end{equation}
On the other hand, by choice of the stopping times  $ \rho  , $
we have 
\begin{equation}\label{him3} \|   R_k \| \le  2 \lambda, 
\end{equation}
Use Doob's martingale projection to generate analytic functions.
Define $$
 \Delta G_{k} 
=  N( R_k ) 
, \quad\quad 
 \Delta B_{k}  
 =     N( S_k )  
, $$
%
where $N$ acts on the last  variable of $ S_k , R_k . $ 
Clearly $ \Delta F_{k} = \Delta G_{k}  + \Delta B_{k} ,
 $ and since Doob's projection preserves analyticity,  $(G_{k}) $ and $ ( B_{k} ) $ form 
Hardy martingales. 
By convexity, the interpretation of Doob's projection  $ N $ as a conditional expectation, 
together with \eqref{him2}, and \eqref{him1} gives 
$$ \bE ( \sum _{k= 1 }^n  \|\Delta B_k\|_X )  \le 
  \bE  ( \sum _{k= 1 }^n  \|S_k\|_X )  \le 8  \bE (  F_n ^* )   \le C \bE  (\|  F_n \|_X   )    . $$
Using once again that   Doob's projection  $ N $ acts as a conditional expectation operator, we get with \eqref{him3}
$$\| \Delta G_{k} \|_X  =    \|N( R_k ) \|_X \le   2 \max_ {m \le k-1} \|F_m ( x ) \|_X  $$
\endproof
 
\subsection{Illustration: Extrapolation of Hardy-Martingale Transforms}
Throughout this section we fix  a Banach space $X$ and  $  \varepsilon =  ( \varepsilon _m ) 
\in \{ -1, 1 \} ^\bN . $
We 
  define  the  operators 
\begin{equation}\label{linz2}
 T_ \varepsilon ( F ) = \sum  \varepsilon _m \Delta F _m , \quad \quad   T_ \varepsilon ( F )_k  = \sum_{m = 1 }^k  \varepsilon _m \Delta F _m . 
 \end{equation}
initially for 
finite, $X$ valued Hardy martingales $ F = ( F_k ) _{k = 1} ^n . $

Illustrating how the Davis decomposition for vector valued 
Hardy martingales may be applied, we combine it with  
an extrapolation method for {\it previsible} martingales (Maurey \cite{m2}.)
Thus Theorem \ref{dc1} yields  Garling's  \cite{gar2}  extrapolation theorem.  
\begin{theor}\label{dc2}
If  there exists $ A_2 > 0 $ such that for any 
square integrable, $X$ valued Hardy martingale $ Z = ( Z_k ) $  
$$\bE ( \|  T_ \varepsilon ( Z ) _k \| _X^2  ) \le A_2 ^2 \bE ( \| Z  _k \| _X^2  ), \quad \quad k \in \bN ,$$
then there exists $A_1= A_1 ( A_2) $ such that  for any integrable $X$ valued Hardy martingale $ F = ( F_k ) $  
$$\bE ( \|  T_ \varepsilon ( F) _k \| _X  ) \le  A_1   \bE ( \| F  _k \| _X  ), \quad \quad k \in \bN . $$
\end{theor}
\paragraph{Remarks:}
The proof by  Garling  \cite{gar2} combined  extrapolation for previsible 
martingales (e.g.  Burkholder \cite{bur1}) and used that  Q.  Xu has shown 
that Edgar's  approximation argument ( \cite{edgar1}, \cite{edgar2} )  reduces the problem to 
a special case, called analytic martingales. 
For a recent study of the operators  $  T_ \varepsilon $ we refer to the results in the thesis of Yanqi Qiu \cite{qiu-yanqi-diss, qiu-yanqi}.

We first recall 
Maurey's extrapolation argument \cite{m2}. 
\begin{lemma} (Maurey \cite{m2}.)  \label{dc3}
Assume that there exists $ A_2 > 0 $ so that for any 
square integrable, $X$ valued Hardy martingale $ Z = ( Z_k ) $  
$$\bE ( \|  T_ \varepsilon ( Z ) _k \| _X^2  ) \le A_2 ^2 \bE ( \| Z  _k \| _X^2  ) \quad \quad k \in \bN .$$
Let   $ G = ( G_k ) $  be an   $X$ valued  integrable Hardy martingale. Let 
$ w = ( w_k ) $ be a non negative, increasing and adapted sequence satisfying
\begin{equation}\label{linz4} \max _{ m \le k } \| G_k \| _X   \le w_{k-1} .
 \end{equation} 
Then 
 $$\bE ( \|  T_ \varepsilon ( G) _k \| _X  ) \le 8 A_2   \bE ( w_{k-1}  ) , \quad \quad k \in \bN . $$
\end{lemma}
\proof 
We follow the basic steps of Maurey's argument in  \cite{m2}.
\paragraph{Step 1.} Given   $ G = ( G_k ) $ define the transformed Hardy martingale 
$$ Z_k = \sum_{m = 1 } ^k w_{m-1}^{-1/2} \Delta G_m ,  \quad \quad k \in \bN .  $$
\paragraph{Step 2.} We infer from Maurey \cite{m2} that with \eqref{linz4}, the transformed Hardy martingale  $ Z = ( Z_k ) $ 
satisfies the pointwise estimates  
\begin{equation}\label{him4} \|  Z_k \| _X \le  2  w_{k-1}^{1/2} , \end{equation} 
and 
 \begin{equation}\label{him5} \|   T_ \varepsilon ( G )_k  \| _X \le  2 (\max_{ m \le k }\|     T_ \varepsilon ( Z )_m \|_X )    w_{k-1}^{1/2}  .\end{equation} 
\paragraph{Step 3.} By  \eqref{him5},   the Cauchy Schwarz inequality and Doob's maximal theorem we obtain 
$$  \bE (    \|   T_ \varepsilon ( G )_k  \| _X  ) \le  2  ( \bE  \|  \max_{ m \le k }   T_ \varepsilon ( Z )_m \|_X ^2 ) ^{1/2}  ( \bE (   w_{k-1} ) )  ^{1/2} \le 4        ( \bE  \|     T_ \varepsilon ( Z )_k \|_X ^2 ) ^{1/2}  ( \bE (   w_{k-1} ) )  ^{1/2} . $$ 
Next, by  the  hypothesis on $ T_ \varepsilon $ and the pointwise bound \eqref{him4},  we get 
$$  \bE ( \|     T_ \varepsilon ( Z )_k \|_X ^2 ) \le A_2^2  \bE ( \|     Z _k \|_X ^2   ) \le 
  4  A_2^2       \bE (   w_{k-1} ) . $$ 
Summing up we have 
$$  \bE (    \|   T_ \varepsilon ( G )_k  \| _X  ) \le   8  A_2    \bE (   w_{k-1} ) . $$ 

\endproof

\paragraph{Proof of Theorem \ref{dc2}.}
With Theorem \ref{dc1} decompose the Hardy martingale as $F = G+ H . $
We use Lemma \ref{dc3} to estimate  $ T_ \varepsilon ( G )$
and the triangle inequality for  $ T_ \varepsilon ( B )$.
\paragraph{Step 1.} Apply 
Theorem \ref{dc1} to the X valued Hardy martingale $ F = ( F_k) $ and obtain  the splitting as $ F = G + H $ 
such that 
\begin{equation}\label{linz3}
\| \Delta G _ k \| _X \le 2 \max_{ m \le k -1 }   \| F_m \| _X ,
\quad\text{and}\quad 
 \bE ( \sum _{ m = 1 } ^k  \| \Delta B_ m \| _X ) \le C_0  \bE (   \| F_ k  \| _X ) ,  \end{equation}
where again $ G, B $ are   X valued Hardy martingales.
\paragraph{Step 2.} Put 
$$ w_{k-1} =  2 \max_{ m \le k -1 }   \| F_m \| _X  +  \max_{ m \le k -1 }   \| G_m \| _X  . $$ 
By \eqref{linz3},   $   \|  G _ k \| _X   \le w_{k-1 }  .$  Hence  
Lemma \ref{dc3} applies and gives
 \begin{equation}\label{him6}   \bE (    \|   T_ \varepsilon ( G )_k  \| _X  ) \le   8  A_2    \bE (   w_{k-1} ) .\end{equation}  
By \eqref{linz3}, and the maximal function estimates for Hardy martingales  in \eqref{him0} we have 
 \begin{equation}\label{him7}   \bE (   w_{k-1} )  \le C_1   \bE (   \| F_ k  \| _X  ) .\end{equation} 
  \paragraph{Step 3.}
Next we turn to estimating $  T_ \varepsilon ( B )_k  . $  We use  \eqref{linz3} and triangle inequality as follows, 
 \begin{equation}\label{him8}  \bE (    \|   T_ \varepsilon ( B )_k  \| _X  ) \le   \bE ( \sum _{ m = 1 } ^k  \| \Delta B_ m \| _X ) \le C_0  \bE (   \| F_ k  \| _X ) .\end{equation}  
Summing up the estimates \eqref{him6} --  \eqref{him8} we  get 
 $$   \bE (    \|   T_ \varepsilon ( F )_k  \| _X  ) \le  \bE (    \|   T_ \varepsilon ( G )_k  \| _X  )  +  \bE (    \|   T_ \varepsilon ( B )_k  \| _X  ) \le  ( 8  A_2  C_1 + C_0   ) 
 \bE (   \| F_ k  \| _X ) . $$ 
\endproof

\subsection{The Strong Davis Decomposition}
We continue with decomposition theorems. In Theorem \ref{dc4} 
we determine  a  splitting of a 
 vector valued Hardy martingale $F$  as $F = G + B $ that improves apon  the 
classial Davis decomposition of Theorem \ref{dc1}. 
Specifically   the uniform estimates for the 
predictable part  $G  $  are improved.

In addition to Brownian motion and stopping times, the proof below 
 makes use of   Havin's Lemma for which we refer to  A.  Pelczynski \cite{pelc} and  J. Bourgain \cite{b1}.

\begin{theor}\label{dc4} Let $X$ be a Banach  space. 
Any  $X$ valued  Hardy martingale  $F = (F_k) $
can be decomposed into the  sum of $X$ valued   Hardy martingales $ F = G + B $ 
such that
$$ \|\Delta G_k \|_X \le C \|F_{k-1}\|_X , $$
and
$$ \bE ( \sum _{k= 1 }^\infty  \|\Delta B_k \|_X ) \le C  \bE( \|F \|_X) . $$
\end{theor} 
The  splitting of the Hardy martingale $F$   
is done separately  for each  martingale difference. Here the proof relies on 
the following decomposition  theorem for vector valued analytic functions.
\begin{theor}\label{dc5} 
For any $ h \in H^1_0  (\bT , X ) $ and $  z \in X $ there exists $ g \in 
  H^\infty _0 (\bT , X ) $  so that 
$$ \| g ( \zeta ) \|_X \le C_0 \| z \| _X , \quad \quad \zeta \in \bT  $$
and 
 \begin{equation}\label{him9a}
 \| z \| _X  + \frac{1}{8} \int_\bT  \| h - g \|_X  dm \le \int_\bT \| z + h  \|_X  dm .  \end{equation} 
The constant satisfies  $ C_0 \le 24 . $  
\end{theor}
\proof The proof begins with the definition of  $  g \in 
  H^\infty _0 (\bT , X ) .$ Thereafter we successively collect the lower 
estimates for the right hand side of \eqref{him9a}.   
\paragraph{Step 1.}
Determine  $  g \in 
  H^\infty _0 (\bT , X ) $  
by putting
$$ \rho = \inf \{ t < \t : \| h( z_t ) \|_X > C_0  \| z \| _X  \} , \quad \quad g = N( h( z_\rho ) ) , $$ 
where $N$ denotes  Doob's projection operator. Since $ h \in H^1_0  (\bT , X ) $
we have 
$ \bE (  h( z_\rho ) ) = 0 $
By definition of the stopping time $  \rho $  we have 
the uniform estimate $ \|  h( z_\rho ) \| _X \le  C_0  \| z \| _X $ and  we obtain 
$$ \| g ( \zeta ) \|_X \le C_0 \| z \| _X , \quad \quad \zeta \in \bT , $$
because Doob's projection $N$ beeing a conditional expectation, is a contraction 
between  $L^\infty $ spaces.  Finally since $N$ preserves analyticity,   we get  $  g \in 
  H^\infty _0 (\bT , X ) . $

\paragraph{Step 2.} We now turn to proving the integral estimates. The idea is to find a lower estimate for the 
right hand side by integrating it against a suitable testing functions.  
Define $ A = \{ \rho < \t \} .$  The set $A$ is measurable with respect to the stopping time 
$ \s -$algebra  $\cF_\rho . $ Since conditional expectations are $L^1 $ contractions we have
 \begin{equation}\label{him11}
\bE ( \| h ( z _ \t ) 1_A \| _X ) \ge  \bE ( \| h ( z _ \r ) 1_A \| _X )  \ge C_0   \| z \| _X . \end{equation} 
Next define
$$ p =  N ( 1_A ) / 2 . $$
Clearly $ 0 \le p \le 1/2 $ and by the covariance formula we get 
$ \int p dm = \bP ( A ) / 2 , $
and  
$$ \int _\bT \| h \|_X  p dm = \frac12    \bE ( \| h ( z _ \t 1_A \| _X ) . $$  
Combining this with \eqref{him11} triangle inequality gives 
 \begin{equation}\label{him12} \int _\bT \|z+  h \|_X  p dm  \ge 
\int _\bT \| h \|_X  p dm -   \frac12 \bP ( A ) \| z \| _X  \ge 
( \frac12 - \frac{1}{ 2 C_0 } )  \bE ( \| h ( z _ \t ) 1_A \| _X ) . \end{equation}   
\paragraph{Step 3.}
Let $ q \in H^\infty ( \bT ) $ be the outer function given by 
$$ q = \exp ( \ln ( 1 - p ) + i H \ln ( 1-p ) ) . $$
 Since $ h \in H^1_0  (\bT , X ) $ we have 
$ \int_\bT h q dm= 0 . $
Put $ q_2 = \Im q $ and $q_1 = \Re q .  $ Then by inspection $ \int q_2 dm = 0 $ and 
$$    \int _\bT ( z+  h )  q dm  =   z  \int_\bT  q_1 dm  . $$
Below we will verify that for $ C_1 > 3 $ 
 \begin{equation}\label{him13}  \int_\bT  q_1 dm > 1 - C_1\bP ( A ) .  \end{equation}  
Assuming the crucial estimate \eqref{him13}  for the moment we may continue our chain of inequalities  as follows. 
 \begin{equation}\label{him14}   \int _\bT \|z+  h \|_X  |q| dm  \ge  \left\| \int _\bT ( z+  h )  q dm \right\|_X \ge 
 \| z \| _X (  1 - C_1\bP ( A ) ) . \end{equation} 
Finally we observe that  $ p + |q| = 1 $ and  take the sum of \eqref{him12} and  \eqref{him14}
to obtain
 \begin{equation}\label{him14a}   \int _\bT \|z+  h \|_X  dm \ge   \| z \| _X  +
 \left( \frac12 - \frac{1}{2 C_0 } - \frac{C_1}{C_0 } \right)  \bE ( \| h ( z _ \t ) 1_A \| _X ) . 
 \end{equation}   
\paragraph{Step 4.} Here we prove that 
\begin{equation}\label{him15}  \int_\bT  \| h - g \|_X  \le 2  \bE ( \| h ( z _ \t ) 1_A \| _X ) . \end{equation}    
 As $ A = \{ \rho < \t \} ,$ the following indentity holds
\begin{equation}\label{him16} 
   h ( z _ \t ) - h ( z _ \r )     =     ( h ( z _ \t ) - h ( z _ \r ) )1_A .  \end{equation} 
Using  that Doob 's projection contracts $L^1 $ spaces, we derive from \eqref{him16}  that 
$$ \int_\bT  \| h - g \|_X  dm  \le \bE ( \|   h ( z _ \t ) - h ( z _ \r ) \|_X 1_A ) . $$
Next use  the right hand inequality in \eqref{him11}  to get 
$$ \bE ( \|   h ( z _ \t ) - h ( z _ \r ) \|_X 1_A )  \le 2  \bE ( \|   h ( z _ \t )   1_A  \|_X ) .    $$ 

\paragraph{Summing up.}
Choose now $ C_0 \ge 8 C_1 $ so that $( 1/2 - 1/( 2C_0 ) - C_1 / C_0 ) > 1 / 4  $ and merge the 
inequalities \eqref{him14a} -- \eqref{him15}  to obtain 
 $$ 
 \int_\bT \| z + h  \|_X  dm   \ge  \| z \| _X  + \frac{1}{8} \int_\bT  \| h - g \|_X  dm .
$$ 
\endproof
\paragraph{A final  remark.} Here we isolate Bourgain's  idea  \cite{b1} used to  prove that  $ 0 \le p \le 1/2 $ implies that \eqref{him13} holds true. We show 
\begin{equation}\label{him20}  \int_\bT  q_1 dm > 1 -  3 \int_\bT  p dm  . \end{equation} 
Recall that $ q_1 = ( 1-p) \cos ( H( \ln ( 1 -p) ) . $ Use  $ \cos  (x) \ge 1  - x^2 / 2 $ 
to get  the pointwise inequality
 \begin{equation}\label{him21} q_1 = ( 1-p)   - 
\frac12   ( H(  ( 1 -p) ) ) ^2 .  \end{equation}
We thus reduced the $L^1 $ estimate for $ q_1 $ to an $L^2 $ estimate for the Hilbert transform.
Clearly we have 
$$ \int_\bT    ( H( \ln  ( 1 -p) ) ) ^2   dm \le 2 \int_\bT   (\ln ( 1 -p) )  ^2  dm . $$
Now if   $ 0 \le p \le 1/2 $ then $ (\ln ( 1 -p) )  ^2  \le 2 p  , $ and hence 
 \begin{equation}\label{him22} \int_\bT    ( H( \ln  ( 1 -p) ) ) ^2   dm  \le 4  \int_\bT  p dm   . \end{equation} 
Combining now  \eqref{him21} and  \eqref{him22} gives \eqref{him20}.
\endproof

\paragraph{Proof of Theorem \ref{dc4}.}
Let  $k \in \bN  $ and condition on   $\cF_{k-1}  $
by fixing   $x  \in \bT^{k-1} .$
For  $ y \in \bT  $
put 
$$ h(y) = \Delta F_k (x, y )\quad\text{and}\quad
z =  F_{k-1} (x).$$ 
We apply  Theorem  \ref{dc5} to $ h \in  H^1 _0 ( \bT, X )  $  and  obtain 
$ g \in H^\infty_0 ( \bT, X ) $, such that 
 \begin{equation*}
 \|g ( \zeta) \|_X \le C_0 \|z\|_X , \quad \quad \zeta \in \bT  \end{equation*} 
and 
 \begin{equation*}
 \|z\|_X + \frac18 
\int_{\bT}\|h-g\|_X dm
\le \int_{\bT} \|z +
h\|_X dm 
. 
\end{equation*} 
Define the splitting of $  \Delta F_k $ by putting
$$ \Delta G_k(x, y ) = g(y),
\quad\quad\text{and} \quad
\Delta B_k(x , y ) = h(y) - g(y).$$
This gives  
$ \Delta F_k =\Delta G_k + \Delta B_k , $
with
$ \|\Delta G_k \|_X \le C_0\| F_{k-1}\|_X .$
and 
$$ \|F_{k-1}\|_X 
+ \frac18 \bE_{k-1}(\|\Delta B_k\|_X ) \le \bE_{k-1}  (\|F_{k}\|_X) . $$ 
Taking  expectations on both sides and summing the telescoping series 
  gives 
$$ \sum  \bE(\|\Delta B_k\|_X)\le  8 \sup  \bE (\|F_{k}\|_X) . $$
 \endproof

\subsection{Illustration: Vector valued Davis and Garsia Inequality}
Here we show that the strong Davis decomposition  yields  vector valued   Davis and Garsia Inequalities. 
At this point we  need to make an assumption on  the Banach space $X$: Let $ q \ge 2. $
A Banach space $X$ satisfies the hypothesis $\cH(q),$  if for each 
$ M \ge 1 $ there exists $ \d = \d (M) > 0 $ such that for any  $ x \in X $ with $ \| x \| = 1 $ and
$ g \in H^\infty _ 0 ( \bT , X )  $ with $ \| g \| _\infty  \le M , $ 
\begin{equation}\label{6-4-15-5}   \int_\bT \| z + g \|_X  dm  \ge  (  1  + \d \int_\bT  \|  g \|_X ^q dm  )^{1/q} . \end{equation} 
We emphsize that  \eqref{6-4-15-5} is required to hold  only  for uniformly bounded analytic functions $ g ,$  
and that $ \d = \d (M) > 0 $ is allowed to depend on the uniform estimates  $ \| g \| _\infty  \le M . $ 
\begin{theor}\label{6-4-15-4}  Let $ q \ge 2 .$ Let $X$ be a Banach satisfying  $\cH(q)$. There exists $ M > 0 $ 
$\d_q > 0 $ such that 
 for any $ h \in H^1_0  (\bT , X ) $ and $  z \in X $ there exists $ g \in 
  H^\infty _0 (\bT , X ) $  satisfying
\begin{equation}\label{6-4-15-2} \| g ( \zeta ) \|_X \le M  \| z \| _X , \quad \quad \zeta \in \bT , \end{equation} 
and 
 \begin{equation}\label{6-4-15-3}
  \int_\bT \| z + h  \|_X  dm \ge \left(  \| z \| _X ^q  + \d_q \int_\bT  \|  g \|_X^q  dm  \right)^{1/q} + \frac{1}{ 16}  \int_\bT  \| h - g \|_X  dm  .  \end{equation} 
\end{theor}
\proof
Let  $ h \in H^1_0  (\bT , X ) $ and $  z \in X .$ 
By Theoren \ref{dc5}   there exists  $ M \le 24 $ and $ g \in 
  H^\infty _0 (\bT , X ) $  so that 
\begin{equation}\label{6-4-15-10} \| g ( \zeta ) \|_X \le M \| z \| _X , \quad \quad \zeta \in \bT , \end{equation} 
and 
 \begin{equation}\label{6-4-15-11}
\int_\bT \| z + h  \|_X  dm  \ge \| z \| _X  + \frac{1}{8} \int_\bT  \| h - g \|_X  dm .
 \end{equation} 
Next by triangle inequality, 
 \begin{equation}\label{6-4-15-12}
\int_\bT \| z + h  \|_X  dm  \ge    \int_\bT  \| z +g\|_X  dm     - \int_\bT  \| h - g \|_X  dm     . \end{equation} 
 and by hypothesis $\cH (q) $  there exists $ \d = \d ( M ) > 0  $ such that 
\begin{equation}\label{6-4-15-13}
 \int_\bT  \| z +g\|_X  dm \ge 
  (  \| z \| _X ^q  + \d\int_\bT  \|  g \|_X ^ q  dm  )^{1/q} . \end{equation} 
Let $ 0 < \a \le 1 $  form  $ \a   \eqref{6-4-15-11}  +  (1 -\a ) \eqref{6-4-15-12} $ and invoke 
\eqref{6-4-15-13}.  Thus we obtained that  $\int \| z + h  \| $ is larger than the following term,
\begin{equation}\label{6-4-15-14} ( 1-\a )  \| z \| _X + \a   (  \| z \| _X ^q  + \d\int_\bT  \|  g \|_X^ q   dm  )^{1/q} + \frac{( 1- 9\a )}{8}  \int_\bT  \| h - g \|_X  dm     . \end{equation} 
Just by triangle inequality    \eqref{6-4-15-14} is larger than 
$$    (  \| z \| _X ^q  + \d \a^q \int_\bT  \|  g \|_X ^ q  dm  )^{1/q} + \frac{( 1- 9\a )}{8}  \int_\bT  \| h - g \|_X  dm     .  $$
Specifying $ \a = 1/ 18 $ finishes the proof  of \eqref{6-4-15-3}.
\endproof
\begin{theor}\label{6-4-15-26} Let $ q \ge 2 . $ Let $X$ be a Banach satisfying  $\cH(q)$. 
Any  $X$ valued  Hardy martingale  $F = (F_k) $
can be decomposed into the  sum of $X$ valued   Hardy martingales $ F = G + B $ 
such that
$$  \bE (  \sum _{k= 1 }^\infty   \bE_{k-1} (\|\Delta G_k \|_X  
^q ) )^{1/q}    +   \bE ( \sum _{k= 1 }^\infty  \|\Delta B_k \|_X ) \le A_q \bE( \|F \|_X) . $$
\end{theor} 
\proof
Let  $k \in \bN  $ and condition on   $\cF_{k-1}  .$
Fix   $x  \in \bT^{k-1} .$
and  $ y \in \bT  $
and define  
$$ h(y) = \Delta F_k (x, y )\quad\text{and}\quad
z =  F_{k-1} (x).$$ 
We apply  Theorem  \ref{6-4-15-4} to $ h \in  H^1 _0 ( \bT, X )  $ and obtain 
$ g \in H^\infty_0 ( \bT, X ) $, satisfying \eqref{6-4-15-3}.
Substituting back we obtain the decomposing 
$$ \Delta G_k(x, y ) = g(y),
\quad\quad\text{and} \quad
\Delta B_k(x , y ) = h(y) - g(y)$$
such that 
\begin{equation}\label{6-4-15-6}
\bE( \|F_{k-1}\|_X ^q 
+ \d \bE_{k-1}(\|\Delta G_k\|_X ^q)  ) ^{1/q} +   C  \bE(\|\Delta B_k\|_X )  \le   \bE (\|F_{k}\|_X) . 
\end{equation} 
Apply  non-linear telescoping  \cite{b1, pfxm12},  to equation \eqref{6-4-15-6}. This gives 
$$  \bE (  \sum _{k= 1 }^\infty   \bE_{k-1}(\|\Delta G_k \|_X  ^q) )^{1/q}    +   \bE ( \sum _{k= 1 }^\infty  \|\Delta B_k \|_X ) \le A_q 
(\bE( \|F \|_X))^{1/q}  (\bE( \sup_{n \in \bN }\|F _n \|_X))^{1/p},$$
where $ 1/p + 1/q = 1 . $  Invoking \eqref{him0} --the maximal function estimate for Hardy martingales --finishes the proof.
\endproof
\paragraph{Remark:} If in the definition of  $ \cH(q ) $ we had replaced \eqref{6-4-15-5}
by  
\begin{equation}\label{6-4-15-35}   \int_\bT \| z + g \|_X  dm  \ge  (  1  + \d (\int_\bT  \|  g \|_X   dm)^ q  )^{1/q} , \end{equation} 
then the above line of reasoning would have resulted in  the previsible projection estimate 
$$  \bE (  \sum _{k= 1 }^\infty   (\bE_{k-1} \|\Delta G_k \|_X  )^ q    )^{1/q}    +   \bE ( \sum _{k= 1 }^\infty  \|\Delta B_k \|_X ) \le A_q 
\bE( \|F \|_X) . 
$$

\section{Embedding: An Alternative to Decomposing.} 
Our starting point in this section is 
 Maurey's  embedding of  $H^1 ( \bT , X ) $  into 
 Hardy martingales with uniformly small increments. See \cite{m1}.
By iterating Maurey's construction we show that an  arbitrary  Hardy martingale may be considered as  a  subsequence of a Hardy martingale with
the additional property that its increments are dominated by a small,  predictable and increasing  process. 
Our  interest in this result comes from  extrapolation theorems such as Burkholder's \cite{bur1} or Maurey's \cite{m2}. 
As stated above  our Theorem \ref{dc7}  is  probably a weaker version of  the embedding  theorem of Q. Xu, referred to  by  Garling \cite{gar2}. 
Nevertheless with respect to extrapolation Theorem \ref{dc7} allows us to draw similar conclusions.

Let $1/2 >  \e > 0 $ and $w \in \bT ^\bN $ with $ w = ( w_k ) .$ We define inductively  $\vp_1(w) =  \e w_1 ,$
and 
\begin{equation}\label{him24} \vp_{n}(w) =  \vp_{n-1}(w)  + 
\e( 1 - |\vp_{n-1}(w)|)^2 w_n .\end{equation}  
As  proved  by Maurey  \cite{m1} $  \vp = (\vp_n) $
is a  uniformly bounded Hardy martingale whose 
limit  is uniformly distributed over $ \bT , $
that is 
$$
\bP (  \{ w \in \bT^\bN  : \vp ( w ) \in B \} ) =  m ( B ) \quad \quad B \sbe \bT , $$
where $  m ( B ) $ denotes the mormalized Haar measure on $ \bT $.

The following is Maurey's embedding theorem  \cite{m1}.
\begin{theor}\label{dc6} 
For any $ f \in H^1 ( \bT , X) $ 
$$ F_n (w ) = f( \vp _n ( w ) ) , \quad \quad w \in \bT^\bN  $$
defines an   $X$ valued  Hardy martingale for which 
\begin{equation}\label{him25}
 \sup _{n \in \bN  } \bE ( \| F_n \| _X ) = \int_\bT \| f  \|_X   dm 
\end{equation} 
and 
\begin{equation}\label{him26}
 \| \Delta F_n   \| _X \le 2\e \int_\bT \| f  \|_X   dm  .\end{equation}  
\end{theor}
\proof For convenience we sketch Maurey's proof.   
It is straightforward to see that 
$( F_n ) $ is indeed an integrable $X$ valued Hardy martingale and that 
\eqref{him25} holds true. 
 We now turn to  the pointwise estimates \eqref{him26}.
Fix $  w \in \bT^\bN  ,$  and $ n \in \bN . $  Then 
$$ \Delta F_n ( w ) =  f( \vp _n ( w ) ) - f( \vp _{n-1} ( w ) ) . $$ 
Put next $ z =  \vp _n ( w ), u =  \vp _{n-1} ( w ) $ and use the Cauchy integral formula to obtain 
$$ f( z ) - f( u ) = \int_\bT \left\{ \frac { \zeta }{ \zeta - z } -  \frac { \zeta }{ \zeta - u } \right \}
 f( \zeta)  dm(\zeta) .$$
By the triangle inequality we get 
\begin{equation}\label{him30} \|  f( z ) - f( u ) \|_X \le  \frac { | z - u | }{( 1 - |u| ) ( 1 -  | z | }    \int_\bT \|f \|_X  dm  
\end{equation} 
We use the defining recursion \eqref{him24} to see that 
\begin{equation}\label{him27} |\vp_{n}(w) -  \vp_{n-1}(w) | = \e ( 1 - |\vp_{n-1}(w)|)^2\end{equation}   
 and 
\begin{equation}\label{him28}( 1 - |\vp_{n-1}(w)| ) \le ( 1 - \e )^{-1} ( 1 - |\vp_{n}(w)| ) .\end{equation}  
Since we put  $ z =  \vp _n ( w ) $ and $  u =  \vp _{n-1} ( w ) ,$  the relations
 \eqref{him27} and  \eqref{him28} imply that 
\begin{equation}\label{him29} \frac { | z - u | }{( 1 - |u| ) ( 1 -  | z | } \le   2\e .\end{equation}  
Combining the estimates \eqref{him30}  and \eqref{him29}  yields  the following pointwise  bounds for the martingale differences
$$  
 \| \Delta F_n   \| _X \le  2\e  \int_\bT \| f  \|_X   dm  . $$
\endproof

Applying  Maurey's Theorem \ref{dc6} repeatedly we
associate to  
an arbitrary  Hardy martingale   a  subsequence of a Hardy martingale 
{\it with   
small predictable increments}  and almost identical norms. As mentioned above this makes it possible to apply 
standard extrapolation theorems {\it without performing a Davis decomposition.}
In that sense the following embedding  provides an alternative to Hardy-martingale-decomposition.

\begin{theor} \label{dc7} Let $ \eta > 0 $, and $ 1 \le p < \infty . $ For any $X$ valued Hardy martingale 
$ g = ( g_k ) $ there exists an X valued  Hardy martingale $ G = (G_k )  $, an increasing sequence of integers
$$ m( 0 ) < m( 1) < \dots < m(n) < \dots $$
and a non-negative  adapted process  $ ( \b _k ) $ such that 
\begin{equation}\label{2-4-15-1} \bE (\sup_{k \in \bN }  \b_k )  \le \sup_{k \in \bN } \bE (\|g_k\| _X ) , \end{equation}
and so that the following conditions hold:
\begin{enumerate} 
\item Small  and previsible increments,
\begin{equation}\label{2-4-15-2}
 \| \Delta G_k    \| _X  \le \eta  \b _{k -1}, \end{equation}
\item Almost identical $L^p $ norms,
\begin{equation}\label{2-4-15-3}
( 1 - \eta ) \bE ( \| g_k \|_X ^p  ) \le  \bE ( \| G_{m(k)} \|_X ^p  ) \le   ( 1 - \eta ) \bE ( \| g_k \|_X ^p  )  . \end{equation}
and 
\begin{equation}\label{2-4-15-4}
( 1 - \eta ) \bE ( \| \Delta g_k \|_X ^p  ) \le  \bE ( \| G_{m(k)} - G_{m(k-1)}  \|_X ^p  ) \le   ( 1 - \eta ) \bE ( \|  \Delta g_k \|_X ^p  )  . \end{equation}
\end{enumerate}
\end{theor}

\proof
The proof iterates  Maurey's Theorem \ref{dc6}. 
First of all we my  assume that the martingale $ g = ( g_k ) $ is finite,
and that moreover each $ g_k $ is a trigonometric polynomial. 
To keep the notation simple we restrict the presentation 
to the case $ p = 1 . $
\paragraph{Step 1 (Preparation).}
Depending on the  martingale $ ( g_k ) $  we select  $ 0 < \e <\eta  $  
so that 
\begin{equation}\label{him40} \e / \eta =   ( \sup \bE (\| g_k \| _X ) {\Large{/}} (  \sum_{k \in \bN } \bE (  \| \Delta g_k \| _X ) .  
\end{equation}
Since  $ g = ( g_k ) $ is finite we have in fact $ 0 < \e .$ 
Let $ \vp= ( \vp_n ) $ be the 
Hardy martingale \eqref{him24}
defined by $\vp_1(w) =  \e w_1 ,$
and 
$$ \vp_{n}(w) =  \vp_{n-1}(w)  + 
\e( 1 - |\vp_{n-1}(w)|)^2 w_n , \quad\quad w \in \bT^\bN .$$
We shorten the notation and put 
$$ \O = \bT^\bN . $$
\paragraph{Step 2 (Substitution).}
For 
$ k \in \bN $  and 
 $u \in \O ^{k}$ 
 we write  $ u = (  u^{(1) } , \dots ,  u^{(k) })  $
where $   u^{(1) }  \in  \O ,$ $ \dots$,  $ u^{(k) }  \in  \O . $
Define 
 \begin{equation}\label{him41}
\Phi ^{k} : \O ^{k} \to \bT ^{k},\quad\quad \Phi ^{k} ( u) = ( \vp ( u^{(1) })  , \dots ,  \vp ( u^{(k) }) ),
\end{equation}
and form the  linear substitution operator
\begin{equation}\label{him42} 
T : L^1 (  \bT^k ) \to  L^1 (  \O^k ), \quad\quad  
  T f ( u ) = f (  \Phi ^{k} ( u)  ) .
\end{equation}
Clearly, $T$ is a contraction between the   $L^1$ spaces in \eqref{him42}.

Fix $ k \in \bN $, $ v \in   \O^{k-1} ,$ and $ w \in \O .$ Then clearly $ u = (v, w ) \in  
\O^{k}  $ and 
 we have  
$$(Tg_k ) ( v , w ) = g_k (  \Phi ^{k-1} ( v) , \vp ( w )  ) .$$
We view $g_{k-1}$  as a function on $  \bT^k $ that does not depend on the last variable.
Hence we may apply the substitution operator $T$ to $ g_{k-1}$, and  since 
$ \bE_{k-1} ( g_k ) =  g_{k-1} $ we observe  the following commutation relation 
between expectations
\begin{equation}\label{him95} 
 ( Tg_{k-1} ) ( v ) = \bE _{( w) } ( T g_k ( v , w ) ) . 
\end{equation}
\paragraph{Step 3 (An intermediary Hardy martingale).} We fix   $ v \in   \O^{k-1} $ and form the   Hardy martingale with respect to the last variable,
$$ h_m ( w ) = g_k (  \Phi ^{k-1} ( v) ,  \vp _m  ( w )  ) , \quad \quad m \in \bN , w \in \O . $$
Theorem  \ref{dc6} asserts that $ h = ( h_m ) $ is an X valued  Hardy martingale,
and that its  increments are small and  uniformly bounded. 
Specifically, if we put   
$$\a _{k-1} ( v ) =  \bE _{( w) }  ( \| Tg_k ( v , w ) - T g_{k-1}  (v) \| _X ) , $$
then 
$$ \sup_{ w \in \O } \| \Delta  h _n (w)  \| _X \le \e  \a _{k-1} ( v) $$
and 
$$   \sup_{n \in \bN }   \bE _{( w) } \|  h _n \| _X   =      \bE _{( w) }  ( \| Tg_k ( v , w ) \| _X ) , $$
where the integration is over $ w \in \O
$

\paragraph{Step 4 (Bounding the active variables).} Since $g = ( g_k) $ is assumed to be a finite martingale we pick now $ n \in \bN $ so that 
 \begin{equation}\label{him90} g = ( g_k) _{k = 1 } ^ n 
\end{equation} 
We approximate $Tg_k $ by stopping the martingales $ \vp = ( \vp _m ) $
used in the definition of the linear substitutions $T $.  We will  replace the limit  $ \vp $ by 
one of its  approximatants $ \vp _m ,$ thereby reduce  the number of active variables.

For any  $m \in \bN $ define the substitutions 
$$ T _m  f ( v ) = f (  \Phi _m  ^{n} ( v)  ) ,\quad\quad  \Phi_m  ^{n} ( v) = ( \vp_m ( v^{(1) })  , \dots ,  \vp_m  ( v^{(n) }) ),\quad\quad v \in \O ^{n} ,$$
 Since  $ n \in \bN $ and $ \e > 0 $ are fixed, there exists  $ K \in \bN $ so that 
  \begin{equation}\label{him91}
\sup _{k \le n } \bE_{ \O ^n} ( \| T _K  ( g_k)  -    T  ( g_k) \|_X )  \le \e  \sup _{k \le n } 
\bE_{ \O } ( \|  g_k \|_X ) , \end{equation}  
and  \begin{equation}\label{2-4-15-6}
\sup _{k \le n } \bE_{ \O ^n} ( \| T _K  ( \Delta g_k ) -    T  (  \Delta g_k ) \|_X )  \le \e  \sup _{k \le n } 
\bE_{ \O } ( \|  g_k \|_X ) , \end{equation} 
where the integration on the left hand side is with respect to the normalized Haar measure of 
 $ \O ^n  ,$ and on the right hand side we integrate over $ \O =  \bT^\bN . $ 

By construction, for each $ k \le n ,$ the dependence of 
$ T_K ( g_k ) $ is only on the following variables, 
$$ ( v^{(1)}  _1  , \dots ,   v^{(1)} _{K}  , \dots \dots ,  
v^{(k-1)}_1   , \dots ,   v^{(k-1)} _{K} ,  w_1    , \dots ,  w _{K} )  .$$ 
Thus $ T_K ( g_k ) $ is un-ambiguously defined on the torus product 
$$ \bT ^{K k } \sbe \bT ^{K n } . $$
\paragraph{Step 5 (Conclusion).}
Put   $ N = K n , $ where $ n \in \bN $ respectively   $K \in \bN$  are  defined
by  \eqref{him90} respectively \eqref{him91}.
Finally we define the Hardy martingale $ G = ( G_k )_{k = 1 } ^N  :$ 
 Put 
$$ 
\rho  : \O ^n \to \bT ^{Kn } , \quad\quad    \rho   ( v) = 
 ( v^{(1)}  _1  , \dots ,   v^{(1)} _{K }  , \dots  ,  v^{(n)}_1   , \dots ,   v^{(n)} _{K}  ) ,$$
then 
$$G : \bT ^{ N } \to X $$
is defined  without ambiguity, by putting
$$ G( z ) = ( T_K ( g_n ) ( v) , \quad\quad  z =  \rho   ( v) . $$ 
Let  $ m ( k ) = K k $  then  by the commutation relation \eqref{him95} 
$$ G_{  m ( k )} =  \bE _{  m ( k )} ( G ) ( z ) =  ( T_K ( g_k ) ( v) , \quad\quad  z =  \rho   ( v)  , $$
and 
$$   (  G_{  m ( k )}    -  G_{  m ( k -1)} ) ( z )  =  T_K (  g_k  -  g_{k-1} ) ( v) , \quad\quad  z =  \rho   ( v)  . $$
Thus in view of \eqref{him91} and \eqref{2-4-15-6} we verified  \eqref{2-4-15-3} and  \eqref{2-4-15-4}.

Finally we let $ \bE _K $ denote the conditional  expectation projecting onto the 
first $ K $ variables of $ \O = \bT^\bN . $ Let 
$$ \bF_K =  \bE _K \otimes \dots \otimes \bE _K  , $$
be the conditional expectation  on $ \O^n $ 
given by the $ n -$ fold tensor product of $ \bE _K $. 
Theorem \ref{dc6} asserts that for  $ m( k-1) \le j < m(k)  ,  $ we have the pointwise estimate
 \begin{equation}\label{2-4-15-7}
 \| \Delta G_j ( z ) \| _X \le \e  \bF_K ( \a _{k-1}  ( v ) ) ,  \quad\quad  z =  \rho   ( v)  . \end{equation}
Define now the adapted process 
$$  \b _{k-1} ( z )  = ( \e / \eta)  \bF_K ( \a _{k-1}  ( v )  ,  \quad\quad  z =  \rho   ( v)  . $$
Clearly we have 
$
 \bE_{ \O ^n} ( \sup _k  \bF_K ( \a _{k-1}   ) )  \le   \bE_{ \O ^n} ( \sum \|   T  \Delta g_k  \|_X )  , $ 
and   \eqref{him40} --specifying    the relation between  $ \e $ and $  \eta > 0 $--gives  
 \begin{equation}\label{2-4-15-8} \bE (\sup_{k \in \bN }  \b_k )  \le \sup_{k \in \bN } \bE (\|g_k\| _X ) .\end{equation}
Thus  \eqref{2-4-15-7}  translates to   \eqref{2-4-15-2} and  \eqref{2-4-15-8}  gives  \eqref{2-4-15-1}. 

\endproof

\bibliographystyle{abbrv}
\bibliography{hardymartingales}
Department of Mathematics\\
J. Kepler Universit\"at Linz\\
A-4040 Linz\\
paul.mueller@jku.at

\end{document}